\documentclass[a4paper,12pt]{article}
\pdfoutput=1
   \usepackage[top=2.5cm,bottom=2.5cm,left=2.5cm,right=2.5cm]{geometry}
\usepackage{cite,amsmath,amssymb}
\usepackage[margin=1cm,%
            font=small,%
            format=hang,%
            labelsep=period,%
             labelfont=bf]{caption}
\usepackage{graphicx}
\usepackage{caption}
\usepackage{float}
\usepackage{amsmath}
\usepackage{cases}
\usepackage{mathrsfs}
\usepackage{amsfonts}
\usepackage{color}
\usepackage{setspace}
\usepackage{cite}
\usepackage{amsmath}
\usepackage{array}
\usepackage{booktabs}
\usepackage{graphicx}
\usepackage{subfigure}
\usepackage{CJK}
\usepackage{caption}
\usepackage{float}
\usepackage{amsmath}
\usepackage{cases}
\usepackage{mathrsfs}
\usepackage{amsfonts}
\usepackage{times}
\usepackage{mathptmx}
\usepackage{cite}
\usepackage{setspace}
\usepackage{times}
\usepackage{mathptmx}
\usepackage{geometry}

\renewcommand{\figurename}

\renewcommand{\figurename}

\makeatletter
\def\thanks#1{\protected@xdef\@thanks{\@thanks
        \protect\footnotetext{#1}}}
\makeatother

\begin{document}
\title{\textbf{Extremal trees with respect to spectral radius of restrictedly weighted adjacency matrices}
\vspace*{0.3cm}}
\author{ \bf  Ruiling Zheng,\, Xiaxia Guan \, Xian$^{^{\textbf{,}}}$an Jin\\
\small\em School of Mathematical Sciences, Xiamen University, Xiamen  361005, P. R. China\\ \small\em E-mail:
rlzheng2017@163.com, gxx0544@126.com, xajin@xmu.edu.cn.\vspace*{0.3cm}}
 \date{}
\maketitle\thispagestyle{empty}
\vspace*{-1cm}

\begin{abstract}
For a graph $G=(V,E)$ and $v_{i}\in V$, denote by $d_{i}$ the degree of vertex $v_{i}$. Let $f(x, y)>0$ be a real symmetric function in $x$ and $y$. The weighted adjacency matrix $A_{f}(G)$ of a graph $G$ is a square matrix, where the $(i,j)$-entry is equal to $\displaystyle f(d_{i}, d_{j})$ if the vertices $v_{i}$ and $v_{j}$ are adjacent and 0 otherwise. Li and Wang \cite{U9} tried to unify methods to study spectral radius of weighted adjacency matrices of graphs weighted by various topological indices. If $\displaystyle f'_{x}(x, y)\geq0$ and $\displaystyle f''_{x}(x, y)\geq0$, then $\displaystyle f(x, y)$ is said to be increasing and convex in variable $x$, respectively. They obtained the tree with the largest spectral radius of $A_{f}(G)$ is a star or a double star when $f(x, y)$ is increasing and convex in variable $x$. In this paper, we add the following restriction:
$f(x_{1},y_{1})\geq f(x_{2},y_{2})$ if $x_{1}+y_{1}=x_{2}+y_{2}$ and $\mid x_{1}-y_{1}\mid>\mid x_{2}-y_{2}\mid$
 and call $A_f(G)$ the restrictedly weighted adjacency matrix of $G$. The restrictedly weighted adjacency matrix contains weighted adjacency matrices weighted by first Zagreb index, first hyper-Zagreb index, general sum-connectivity index, forgotten index, Somber index, $p$-Sombor index and so on. We obtain the extremal trees with the smallest and the largest spectral radius of $A_{f}(G)$. Our results push ahead Li and Wang's research on unified approaches.

\medskip

\end{abstract}

\vspace*{0.35cm}
\baselineskip=0.30in

\section{Introduction}

\noindent

Let $G=(V(G),E(G))$ be a finite, undirected, simple and connected graph with vertex set $V(G)=\{v_{1},v_{2},\ldots,v_{n}\}$ and edge set $E(G)$. An edge $e\in E(G)$ with end vertices $v_{i}$ and $v_{j}$ is usually denoted by $v_{i}v_{j}$. For $i=1,2,\ldots,n$, we denote by $d_{i}$ the degree of the vertex $v_{i}$ in $G$, $\Delta(G)$ the maximum degree of $G$, $N(v_{i})$ the set of neighbours of vertex $v_{i}$ in $G$ and $N[v_{i}]=N(v_{i})\cup\{v_{i}\}$. As usual, let $P_{n}$, $S_{n}$ and $C_{n}$ be the path, star and cycle of order $n\geq3$, $S_{d,n-d}$ be the double star of order $n\geq4$ with two centers $v_{1}$, $v_{2}$ such that $d_{1}=d$ and $d_{2}=n-d$ where $\displaystyle 2\leq d\leq \lfloor\frac{n}{2}\rfloor$.

Let $\displaystyle \lambda_{i}(M) (i=1,2,\ldots,n)$ be the eigenvalues
of a complex matrix $M$. Then the spectral radius of $M$ is
$\displaystyle \rho(M)=\mathop{\text{max}}\limits\mid\lambda_{i}(M)\mid$
. By Perron-Frobenius theorem, when $M$ is an $n\times n$ nonnegative and irreducible matrix, its spectral radius $\rho(M)$ will be a simple eigenvalue of $M$ and there exists a positive real vector ${\bf{x}}=(x_{1}, x_{2},\ldots, x_{n})^{\intercal}$ such that $M{\bf{x}}=\rho(M){\bf{x}}$. The positive eigenvector is called the principal eigenvector of $M$.

In molecular graph theory, the topological indices of molecular graphs are used to reflect chemical properties of chemical molecules. There are many topological indices and among them there is a family of degree-based indices.
The degree-based index $TI_{f}(G)$ of $G$ with positive symmetric function $f(x,y)$ is defined as
\begin{center}
$ \displaystyle TI_{f}(G)=\sum\limits_{v_{i}v_{j}\in E(G)}f(d_{i}, d_{j}).$
\end{center}
There are many degree-based indices such as Randi\'{c} index \cite{1}, Atom-Bond Connectivity index \cite{8}, Arithmetic-Geometric index \cite{15} and Sombor index \cite{U4} and so on.

Each index maps a molecular graph into a single number. Li \cite{U7} proposed that if we use a matrix
to represent the structure of a molecular graph with weights separately on its pairs of adjacent vertices, it will keep more structural information of the graph. For example, the Randi\'{c} matrix \cite{555,666}, Atom-Bond Connectivity matrix \cite{10}, Arithmetic-Geometric matrix \cite{16} and Sombor matrix \cite{U6} were considered separately. Based on these examples, Das et al. \cite{U8} proposed the weighted adjacency matrix $A_{f}(G)$, and it is defined as

$$A_{f}(G)(i,j)=\left\{
\begin{aligned}
&f(d_{i}, d_{j}), \ \ \ \  v_{i}v_{j}\in E(G);\\
&0,\ \ \ \ \ \ \ \ \ \ \ \ \ otherwise.
\end{aligned}
\right.
$$

Since $G$ is a connected graph, the weighted adjacency matrix $A_{f}(G)$
is an $n\times n$ nonnegative and irreducible symmetric matrix. Thus $\rho(A_{f}(G))$ is exactly the largest eigenvalue of $A_f(G)$ and it has a positive
eigenvector ${\bf{x}}=(x_{1}, x_{2},\ldots, x_{n})^{\intercal}$.  Throughout this paper, we choose ${\bf{x}}$ such that $\|{\bf{x}}\|_{2}=1$ and $x_{i}$ corresponds to vertex $v_{i}$, and we call the unique unit positive vector ${\bf{x}}$ principal eigenvector of $G$.

In the recent years, it is a trend to develop unified methods to deal with extremal problems for such degree-based indices and weighted adjacency matrices, \cite{U13,U16,U14,U15,U10,U11,U12,U9} and the survey \cite{10161}.
If $\displaystyle f'_{x}(x, y)\geq0$ and $\displaystyle f''_{x}(x, y)\geq0$, $\displaystyle f(x, y)$ is said to be increasing and convex in variable $x$, respectively. Hu et al. \cite{U14} proposed that a function $\displaystyle f(x, y)$ is said to have the property $P$ if for any $x_{1}+y_{1}=x_{2}+y_{2}$ and $\mid x_{1}-y_{1}\mid>\mid x_{2}-y_{2}\mid$, $f(x_{1},y_{1})>f(x_{2},y_{2})$. Among all simple graphs with order $n\geq2$ and size $\displaystyle 1\leq m\leq \frac{n(n-1)}{2}$, they obtained the connected graphs with a vertex $v_{i}$ such that $d_{i}=n-1$ are the maximal graphs for $TI_{f}(G)$ where $f(x,y)$ is increasing and convex in variable $x$, and has the property $P$. Li and Wang \cite{U9} firstly tried to find unified methods to study spectral radius of weighted adjacency matrices of graphs weighted by some topological indices. They obtained the tree with the largest spectral radius of $A_{f}(G)$ is $S_{n}$ or double star $S_{d,n-d}$, when $f(x, y)$ is increasing and convex in variable $x$.

The idea of this paper is also to go along the unified methods. If the real symmetric function $\displaystyle f(x, y)>0$ is increasing and convex in variable $x$ and satisfies that $f(x_{1},y_{1})\geq f(x_{2},y_{2})$ when $\mid x_{1}-y_{1}\mid>\mid x_{2}-y_{2}\mid$ and $x_{1}+y_{1}=x_{2}+y_{2}$, we call $A_{f}(G)$ the restrictedly weighted adjacency matrix of $G$. In this paper, we consider the extremal trees with respect to spectral radius of restrictedly weighted adjacency matrices. We obtain among all trees of order $n$, $P_{n}$ is the unique tree with the smallest spectral radius and $S_{n}$ is the unique trees with the largest spectral radius.

\vskip0.3cm

The paper is organized as follows. In sections 2 and 3, we provide some known results and two new lemmas on edge moving. In sections 4, we deal with trees. In the last section, we summarize our work and discuss the difficulties for further study.

\section{Known results}
\noindent

In this section, we provide the knowledge of matrix theory on nonnegative matrices and some results on spectra of graphs used in the subsequent sections.

\noindent  {\bf{Theorem 2.1}} \cite{19}
Let $A$ and $B$ be both $n\times n$ nonnegative symmetric matrices. Then $\displaystyle \rho(A+B)\geq\rho(A)$. Furthermore, if $A$ is irreducible and $B\neq0$, then $\displaystyle \rho(A+B)>\rho(A)$.

\noindent  {\bf{Theorem 2.2}} \cite{U1}
Let $A$ be an $n\times n$ real symmetric matrix and $B$ be the principal submatrix of $A$. Then $\displaystyle \rho(A)\geq\rho(B)$.

\noindent  {\bf{Theorem 2.3}} \cite{U1}
Let $A$ be an $n\times n$ nonnegative and symmetric matrix. Then $\rho(A)\geq{\bf{x}}^{\top}A{\bf{x}}$ for any unit vector ${\bf{x}}$, with equality holds if and only if $\displaystyle A{\bf{x}}=\rho(A){\bf{x}}$.

\noindent {\bf{Definition 2.1.}} Let $A$ be an $n\times n$ real matrix whose rows and columns are indexed by $X=\{1,2,...,n\}$. We partition $X$ into $\{X_{1},X_{2},...,X_{k}\}$ in order and rewrite $A$ according to $\{X_{1},X_{2},...,X_{k}\}$ as follows:
\begin{equation*}
  A=\begin{pmatrix}
    A_{1,1} & \cdots & A_{1,k}\\
    \vdots & \ddots & \vdots \\
 A_{k,1} & \cdots & A_{k,k}\\
\end{pmatrix},
\end{equation*}

\noindent where $A_{i,j}$ is the block of $A$ formed by rows in $X_{i}$ and the columns in $X_{j}$. Let $b_{i,j}$ denote the average row sum of $A_{i,j}$. Then the matrix $B=[b_{i,j}]$ will be called the \textbf{quotient matrix} of the partition of $A$. In particular, the partition is called an \textbf{equitable partition} when the row sum of each block $A_{i,j}$ is constant.

\noindent  {\bf{Theorem 2.4}} \cite{U1}
Let $A\geq 0$ be an irreducible matrix, $B$ be the quotient matrix of an equitable partition of $A$. Then $\displaystyle \rho(A)=\rho(B)$.

Now we recall two extremal results on spectral radius of adjacency matrix of trees. The spectral radius of the adjacency matrix $A(G)=(a_{i,j})$ of $G$ is referred as the spectral radius of $G$ and it is denoted as $\displaystyle \rho(G)$.

\noindent  {\bf{Theorem 2.5}} \cite{T1}
Let $T$ be a tree of order $n\geq3$. Then
\begin{center}
$\displaystyle 2\cos\frac{\pi}{n+1}=\rho(P_{n})\leq\rho(T)\leq\rho(S_{n})=\sqrt{n-1},$
\end{center}
\noindent the lower (upper) bound is attained if and only if $G\cong P_{n}$ ($G\cong S_{n}$).

For the weighted adjacency matrix, in 2021, Li and Wang obtained:

\noindent {\bf{Theorem 2.6}} \cite{U9}
Let $\displaystyle f(x,y)>0$ be a symmetric real function, increasing and convex in variable $x$. Then the tree on $n$ vertices with the largest spectral radius of $\rho(A_{f}(T))$ is $S_{n}$ or a double star $S_{d,n-d}$ for some $\displaystyle 2\leq d\leq \lfloor\frac{n}{2}\rfloor$.

\section{Two lemmas on edge moving}
\noindent

From now on, $A_{f}(G)$ is always assumed to be the restrictedly weighted adjacency matrix of a graph $G$. New lemmas in this section are both mainly based on Theorem 2.3 and the following lemma is similar to and a modification of Lemma 2.1 in \cite{U9}.

Kelmans \cite{U3} introduced a simple local operation of a graph to describe the relation between edge moving and spectral radius as follows.
\vskip0.2cm
\noindent $\textbf{Definition 3.1.}$ Let $v_{1},v_{2}$ be two vertices of the graph $G$. And we denote: $N_{1}=N(v_{1})-N[v_{2}]$, $N_{2}=N(v_{2})-N[v_{1}]$, $N_{3}=N(v_{1})\bigcap N(v_{2})$. We use the Kelmans operation of $G$ as follows: Replace the edge $v_{1}v_{w}$ by a new edge $v_{2}v_{w}$ for all vertices $v_{w}\in N_{1}$ (as shown in Fig. 1). In general, we will denote the obtained graph by $G'$.

\begin{figure}[H]
  \centering
    \includegraphics[width=8cm]{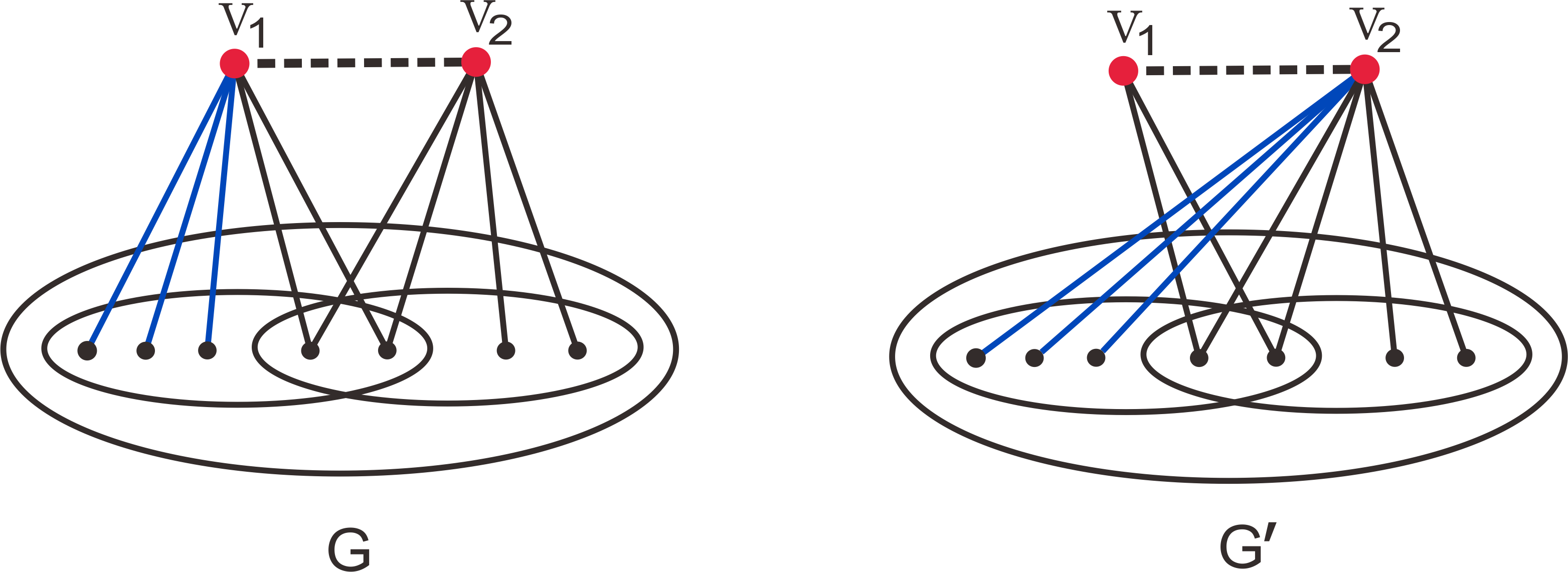}
  \caption*{\small \textbf{Fig. 1.} The Kelmans operation}
\end{figure}

\noindent $\textbf{Remark 3.1.}$ (i) The vertices $v_{1}$ and $v_{2}$ may be adjacent in $G$.

1. Two graphs obtained from moving edges from $N_{1}$ to $N_{2}$ and from $N_{2}$ to $N_{1}$ are isomorphic.

2. $G\cong G'$ if and only if $N_{1}=\emptyset$ or $N_{2}=\emptyset$.

3. $G'$ may have an isolated vertex, thus may be disconnected.

4. Numbers of vertices and edges of $G'$ are the same as those of $G$.

\noindent $\textbf{Lemma 3.1.}$
Let $G$ be a connected graph and $G'$ be the graph after a Kelmans operation on any two vertices $v_{1}$ and $v_{2}$ of $G$ as shown in Fig. 1. If $G\ncong G'$, then $\displaystyle \rho(A_{f}(G))<\rho(A_{f}(G'))$.

\noindent {\textbf{Proof.}} Let $\textbf{x}$ be the principal eigenvector of $G$. By Remark 3.1(ii), without loss of generality, we suppose that $x_{1}\leq x_{2}$.
We denote $n_{1}=\mid N_{1}\mid>0$ and assume the formula
\begin{eqnarray*}
(1)&=&2\sum\limits_{v_{w}\in N_{1}}\big[f(d_{2}+n_{1},d_{w})x_{2}x_{w}-f(d_{1},d_{w})x_{1}x_{w}\big]\\
&+&2\sum\limits_{v_{w}\in N_{2}}\big[f(d_{2}+n_{1},d_{w})-f(d_{2},d_{w})\big]x_{2}x_{w}\\
&+&2\sum\limits_{v_{w}\in N_{3}}\big\{\big[f(d_{2}+n_{1},d_{w})-f(d_{2},d_{w})\big]x_{2}x_{w}\\
&&\quad\quad\quad-\big[f(d_{1},d_{w})-f(d_{1}-n_{1},d_{w})\big]x_{1}x_{w}\big\}.
\end{eqnarray*}
\noindent If $v_{1}v_{2}\notin E(G)$, then
\begin{center}
$\textbf{x}^{\top}A_{f}(G')\textbf{x}-\textbf{x}^{\top}A_{f}(G)\textbf{x}=(1)$.
\end{center}
\noindent If $v_{1}v_{2}\in E(G)$, then
\begin{center}
$\textbf{x}^{\top}A_{f}(G')\textbf{x}-\textbf{x}^{\top}A_{f}(G)\textbf{x}=(1)+ 2\big[f(d_{1}-n_{1},d_{2}+n_{1})-f(d_{1},d_{2})\big]x_{1}x_{2}$.
\end{center}
Because $d_{2}+n_{1}=\mid N_{2}\mid+\mid N_{3}\mid+n_{1}>\mid N_{3}\mid+n_{1}=d_{1}$, $x_{2}\geq x_{1}$ and $f(x, y)$ is increasing and convex in variable $x$, we obtain $(1)\geq0$. Moreover, $f(d_{1}-n_{1},d_{2}+n_{1})-f(d_{2},d_{1})\geq0$. Thus we have
\begin{center}
$\textbf{x}^{\top}A_{f}(G')\textbf{x}-\textbf{x}^{\top}A_{f}(G)\textbf{x}\geq0$.
\end{center}
\noindent According to Theorem 2.3, it follows that
\begin{center}
$\displaystyle \rho(A_{f}(G))=\textbf{x}^{\top}A_{f}(G)\textbf{x}\leq\textbf{x}^{\top}A_{f}(G')\textbf{x}\leq\rho(A_{f}(G'))$.
\end{center}
Let $\displaystyle \rho(A_{f}(G))=\rho(A_{f}(G'))$. Then $\textbf{x}$ is the unit vector such that $\displaystyle A_{f}(G')\textbf{x}=\rho(A_{f}(G'))\textbf{x}$. If $v_{1}v_{2}\in E(G)$, then it can be deduced from

\noindent$\left\{
  \begin{array}{ll}
    \rho(A_{f}(G))x_{2}=\sum\limits_{v_{w}\in N_{2}}f(d_{2},d_{w})x_{w}+\sum\limits_{v_{w}\in N_{3}}f(d_{2},d_{w})x_{w}+f(d_{2},d_{1})x_{1}, \\
    \rho(A_{f}(G'))x_{2}=\sum\limits_{v_{w}\in N_{2}}f(d_{2}+n_{1},d_{w})x_{w}+\sum\limits_{v_{w}\in N_{3}}f(d_{2}+n_{1},d_{w})x_{w}\\
    \ \ \ \ \ \ \ \ \ \ \ \ \ \ \ \ \ +f(d_{2}+n_{1},d_{1}-n_{1})x_{1}+\sum\limits_{v_{w}\in N_{1}}f(d_{2}+n_{1},d_{w})x_{w},
  \end{array}
\right.$

\noindent that $\sum\limits_{v_{w}\in N_{1}}f(d_{2}+n_{1},d_{w})x_{w}=0$, a contradiction. Similarly, if $v_{1}v_{2}\notin E(G)$, we also can deduce $\sum\limits_{v_{w}\in N_{1}}f(d_{2}+n_{1},d_{w})x_{w}=0$, a contradiction.

Thus we have $\displaystyle \rho(A_{f}(G))<\rho(A_{f}(G'))$.
 \hfill$\Box$

According to Lemma 3.1, we have the following corollary.

\noindent \textbf{Corollary 3.1}
Let $G$ be a connected graph which consists of a proper induced subgraph $H$ and a tree $T$ of order $m+1$ $(m\geq 2)$ such that $H$ and $T$ has a unique common vertex $u$ and $T$ is not a star with center $u$. Let $G^{\ast}=H+uv_{1}+\cdots+uv_{m}$, where $v_{1},v_{2},\ldots,v_{m}\in V(T)$ are distinct pendent vertices of $G^{\ast}$. Then $\displaystyle \rho(A_{f}(G))<\rho(A_{f}(G^{\ast}))$.

\noindent {\textbf{Proof.}}
We prove the corollary by induction on the number $p$ of non-pendent vertices different from $u$ in $T$. Because $T$ is not a star with center $u$, we have $p\geq1$.

Assume that $w\in V(T)$ is a non-pendent vertex different from $u$ which is the furthest from $u$ among all non-pendent vertices in $T$. Let $P$ be the unique path in $T$ from $u$ to $w$ and $x$ be the vertex in $P$ which is adjacent to $w$. Moreover, the vertex $x$ may be just $u$. We use a Kelmans operation on the adjacent vertices $w$ and $x$ as follows: Replace the edge $v_{i}w$ by a new edge $v_{i}x$ for all vertices $v_{i}\in N(w)-\{x\}$. Denote the obtained graph by $G_{1}^{\ast}$. Then from Lemma 3.1, we have $\displaystyle \rho(A_{f}(G))<\rho(A_{f}(G_{1}^{\ast}))$. Obviously, the number of non-pendent vertices in the corresponding tree $T_{1}^{\ast}$ in $G_{1}^{\ast}$ different from $u$ is $p-1$. Hence by the induction on $G_{1}^{\ast}$, we have $\displaystyle \rho(A_{f}(G_{1}^{\ast}))\leq\rho(A_{f}(G^{\ast}))$.
Thus we obtain $\displaystyle \rho(A_{f}(G))<\rho(A_{f}(G^{\ast}))$.
 \hfill$\Box$

The graph $F$ is shown in Fig. 2. $v_{1},v_{2}$ are adjacent vertices of $F$, the vertices belong to $N(v_{1})-N[v_{2}]$ and $N(v_{2})-N[v_{1}]$ are all pendent vertices and $v_{z_{1}},\ldots,v_{z_{m}}$ are the common neighbours of $v_{1}$ and $v_{2}$. For convenience, we denote: $N_{1}=N(v_{1})-N[v_{2}]$, $N_{2}=N(v_{2})-N[v_{1}]$, $N_{3}=N(v_{1})\bigcap N(v_{2})$, $n_{i}=\mid N_{i}\mid$ for $i=1,2,3$.

Let $G$ be a connected graph of order $n$ which consists of a proper induced subgraph $H$ and $F$ such that $H$ and $F$ have common vertices $v_{z_{1}},\ldots,v_{z_{m}} (0\leq m\leq n-4)$ and $1\leq n_{1}\leq n_{2}$. Let $\textbf{x}$ be the principal eigenvector of $G$. Assume that the pendent vertices belonging to $N_{1}$ are denoted as $v_{w}$. Due to $A_{f}(G)\textbf{x}=\rho(A_{f}(G))\textbf{x}$. Then the entries corresponding to these pendent vertices in $\textbf{x}$ are equal. The fact also holds for the pendent vertices belonging to $N_{2}$, which are denoted as $v_{y}$. Then we have the following Lemma.

\begin{figure}[H]
\centering
\includegraphics[width=11cm]{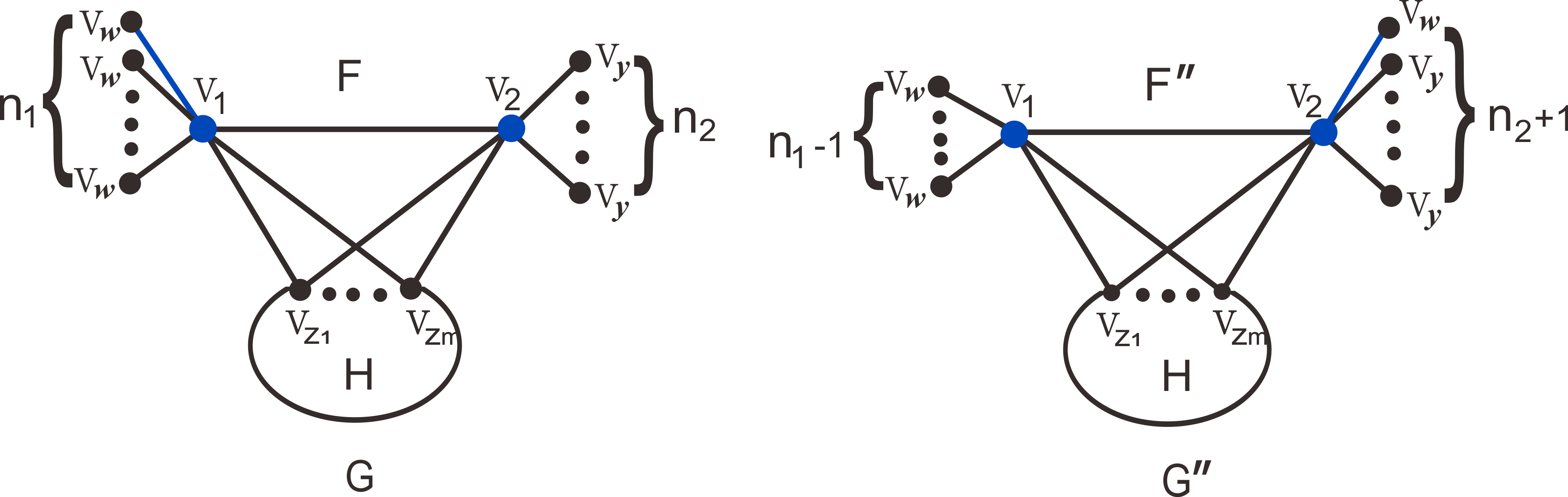}
\caption*{\small \textbf{Fig. 2.} The graphs $F$ and $F''$}
\end{figure}

\noindent \textbf{Lemma 3.2}
 Let $G$ be a connected graph of order $n$ which consists of a proper induced subgraph $H$ and $F$ such that $H$ and $F$ have common vertices $v_{z_{1}},\ldots,v_{z_{m}} (0\leq m\leq n-4)$ and $1\leq n_{1}\leq n_{2}$. Replace $F$ by $F''$ (as shown in Fig. 2) to obtain a new connected graph $G''$.
Then we have $\displaystyle \rho(A_{f}(G))<\rho(A_{f}(G''))$.

\noindent {\textbf{Proof.}}
Due to $A_{f}(G)\textbf{x}=\rho(A_{f}(G))\textbf{x}$, we get

\noindent $\displaystyle \rho(A_{f}(G))x_{1}=\sum\limits_{v_{w}\in N_{1}}f(1,n_{1}+n_{3}+1)x_{w}+\sum\limits_{v_{z_{i}}\in N_{3}}f(d_{z_{i}},n_{1}+n_{3}+1)x_{z_{i}}+f(n_{1}+n_{3}+1,n_{2}+n_{3}+1)x_{2},$

\noindent $\displaystyle \rho(A_{f}(G))x_{2}=\sum\limits_{v_{y}\in N_{2}}f(1,n_{2}+n_{3}+1)x_{y}+\sum\limits_{v_{z_{i}}\in N_{3}}f(d_{z_{i}},n_{2}+n_{3}+1)x_{z_{i}}+f(n_{1}+n_{3}+1,n_{2}+n_{3}+1)x_{1}.$

Hence if $x_{1}>x_{2}$, then $x_{w}>x_{y}$. Meanwhile,
\begin{center}
\noindent $\displaystyle \rho(A_{f}(G))x_{w}=f(1,n_{1}+n_{3}+1)x_{1}$
\end{center}
and
\begin{center}
$\displaystyle \rho(A_{f}(G))x_{y}=f(1,n_{2}+n_{3}+1)x_{2}.$
\end{center}

We obtain if $x_{1}\leq x_{2}$, then $x_{w}\leq x_{y}$.

Firstly, we prove $x_{1}\leq x_{2}$ by contradiction. Assume that $\displaystyle x_{1}>x_{2}$, then we have $\displaystyle n_{2}>n_{1}$. Otherwise, by $A_{f}(G)\textbf{x}=\rho(A_{f}(G))\textbf{x}$, we have $\displaystyle x_{1}=x_{2}$. Replace edge $v_{2}v_{y}$ by a new edge $v_{1}v_{y}$ for $n_{2}-n_{1}$ pendent vertices $v_{y}\in N_{2}$. Denote the obtained graph by $G^{\ast\ast}$. Then we have $G^{\ast\ast}\cong G$.

Since $\displaystyle x_{1}>x_{2}$, $\displaystyle n_{2}>n_{1}$, $x_{w}>x_{y}$ and $\displaystyle f(x,y)>0$ is increasing in variable $x$, we get
\begin{center}
$\big[f(1,n_{2}+n_{3}+1)-f(1,n_{1}+n_{3}+1)\big]x_{1}x_{w}$

$\quad +\big[f(1,n_{1}+n_{3}+1)-f(1,n_{2}+n_{3}+1)\big]x_{2}x_{y}\geq0$,
\end{center}
\begin{center}
 $\displaystyle (n_{2}-n_{1})\big[f(1,n_{2}+n_{3}+1)x_{1}x_{y}-f(1,n_{2}+n_{3}+1)x_{2}x_{y}\big]>0$
\end{center}
\noindent and
\begin{center}
$\sum\limits_{v_{z_{i}}\in N_{3}}\big\{\big[f(d_{z_{i}},n_{2}+n_{3}+1)-f(d_{z_{i}},n_{1}+n_{3}+1)\big]x_{1}x_{z_{i}}$

$\quad  \quad  \quad  \quad  \quad +\big[f(d_{z_{i}},n_{1}+n_{3}+1)-f(d_{z_{i}},n_{2}+n_{3}+1)\big]x_{2}x_{z_{i}}\big\}\geq0$.
\end{center}
\noindent Hence, we have
$$\begin{aligned}
&\quad \textbf{x}^{\top}A_{f}(G^{\ast\ast})\textbf{x}-\textbf{x}^{\top}A_{f}(G)\textbf{x}=\\
&2n_{1}\big[f(1,n_{2}+n_{3}+1)-f(1,n_{1}+n_{3}+1)\big]x_{1}x_{w}\\
+&2n_{1}\big[f(1,n_{1}+n_{3}+1)-f(1,n_{2}+n_{3}+1)\big]x_{2}x_{y}\\
+&2(n_{2}-n_{1})\big[f(1,n_{2}+n_{3}+1)x_{1}x_{y}-f(1,n_{2}+n_{3}+1)x_{2}x_{y}\big]\\
+&2\sum\limits_{v_{z_{i}}\in N_{3}}\big\{\big[f(d_{z_{i}},n_{2}+n_{3}+1)-f(d_{z_{i}},n_{1}+n_{3}+1)\big]x_{1}x_{z_{i}}\\
+&\big[f(d_{z_{i}},n_{1}+n_{3}+1)-f(d_{z_{i}},n_{2}+n_{3}+1)\big]x_{2}x_{z_{i}}\big\}>0.
\end{aligned}$$
This contradicts the fact that $G^{\ast\ast}\cong G$.

Thus $x_{1}\leq x_{2}$ and $x_{w}\leq x_{y}$. Next, we show $\displaystyle \rho(A_{f}(G))<\rho(A_{f}(G''))$.

For $1\leq n_{1}\leq n_{2}$, $x_{1}\leq x_{2}$, $x_{w}\leq x_{y}$ and $\displaystyle f(x,y)>0$ is increasing and convex in variable $x$, we get
\begin{center}
$\displaystyle (n_{1}-1)\big[f(1,n_{1}+n_{3})-f(1,n_{1}+n_{3}+1)\big]x_{1}x_{w}$

$+n_{2}\big[f(1,n_{2}+n_{3}+2)-f(1,n_{2}+n_{3}+1)\big]x_{2}x_{y}\geq0$,
\end{center}

\begin{center}
$\displaystyle f(1,n_{2}+n_{3}+2)x_{2}x_{w}-f(1,n_{1}+n_{3}+1)x_{1}x_{w}\geq0$
\end{center}
\noindent and
\begin{center}
$\displaystyle \big[f(d_{z_{i}},n_{1}+n_{3})-f(d_{z_{i}},n_{1}+n_{3}+1)\big]x_{1}x_{z_{i}}$

$\quad \quad \quad +\big[f(d_{z_{i}},n_{2}+n_{3}+2)-f(d_{z_{i}},n_{2}+n_{3}+1)\big]x_{2}x_{z_{i}}\geq0$.
\end{center}

\noindent Moreover,
\begin{center}
$\displaystyle f(n_{1}+n_{3},n_{2}+n_{3}+2)-f(n_{1}+n_{3}+1,n_{2}+n_{3}+1)\geq0.$
\end{center}
Thus we obtain
$$\begin{aligned}
&\quad \textbf{x}^{\top}A_{f}(G'')\textbf{x}-\textbf{x}^{\top}A_{f}(G)\textbf{x}=\\
&2(n_{1}-1)\big[f(1,n_{1}+n_{3})-f(1,n_{1}+n_{3}+1)\big]x_{1}x_{w}\\
+&2n_{2}\big[f(1,n_{2}+n_{3}+2)-f(1,n_{2}+n_{3}+1)\big]x_{2}x_{y}\\
+&2\big[f(1,n_{2}+n_{3}+2)x_{2}x_{w}-f(1,n_{1}+n_{3}+1)x_{1}x_{w}\big]\\
+&2\sum\limits_{v_{z_{i}}\in N_{3}}\big\{\big[f(d_{z_{i}},n_{1}+n_{3})-f(d_{z_{i}},n_{1}+n_{3}+1)\big]x_{1}x_{z_{i}}\\
+&\big[f(d_{z_{i}},n_{2}+n_{3}+2)-f(d_{z_{i}},n_{2}+n_{3}+1)\big]x_{2}x_{z_{i}}\big\}\\
+&2\big[f(n_{1}+n_{3},n_{2}+n_{3}+2)-f(n_{1}+n_{3}+1,n_{2}+n_{3}+1)\big]x_{1}x_{2}\geq0.
\end{aligned}$$
According to Theorem 2.3, it follows that $\displaystyle \rho(A_{f}(G))\leq \rho(A_{f}(G''))$. If the equality holds, then $\textbf{x}$ is also the principal eigenvector of $G''$. It can be deduced from

\noindent $\left\{
  \begin{array}{ll}
    \rho(A_{f}(G)){x}_{2}=n_{2}f(1,n_{2}+n_{3}+1)x_{y}+f(n_{1}+n_{3}+1,n_{2}+n_{3}+1)x_{1}\\
    \ \ \ \ \ \ \ \ \ \ \ \ \ \ \ \ \  \ \ +\sum\limits_{v_{z_{i}}\in N_{3}}f(d_{z_{i}},n_{2}+n_{3}+1)x_{z_{i}}, \\
    \rho(A_{f}(G'')){x}_{2}=n_{2}f(1,n_{2}+n_{3}+2)x_{y}+f(n_{1}+n_{3},n_{2}+n_{3}+2)x_{1}\\
    \ \ \ \ \ \ \ \ \ \ \ \ \ \ \ \ \  \ \ +\sum\limits_{v_{z_{i}}\in N_{3}}f(d_{z_{i}},n_{2}+n_{3}+2)x_{z_{i}}+f(1,n_{2}+n_{3}+2)x_{w},
  \end{array}
\right.$
 
\noindent that $f(1,n_{2}+n_{3}+2)x_{w}=0$. This leads to a contradiction. Thus $\displaystyle \rho(A_{f}(G))<\rho(A_{f}(G''))$.
 \hfill$\Box$

 \noindent $\textbf{Remark 3.2.}$ The graph $H$ in Lemma 3.2 can be taken to be null.

 According to Lemma 3.2, we have the following corollary.

\noindent \textbf{Corollary 3.2}
Among all double stars of order $n\geq3$, we have $\rho(A_{f}(S_{\lfloor\frac{n}{2}\rfloor,n-\lfloor\frac{n}{2}\rfloor}))<\ldots<\rho(A_{f}(S_{3,n-3}))<\rho(A_{f}(S_{2,n-2}))<\rho(A_{f}(S_{n}))$.

\section{Extremal trees}
\noindent

In this section, we study extremal trees with respect to spectral radius of restrictedly weighted adjacency matrices, we obtain:

\noindent \textbf{Theorem 4.1}
 Among all trees of order $n$, $P_{n}$ is the unique tree with the smallest spectral radius of $A_{f}(G)$ and $S_{n}$ is the unique tree with the largest spectral radius of $A_{f}(G)$.

\noindent {\textbf{Proof.}}
According to Theorem 2.6 and Lemma 3.1, it is easy to show $S_{n}$ is the unique tree with the largest spectral radius of $A_{f}(G)$.

Now we prove $P_{n}$ is the unique tree with the smallest spectral radius of $A_{f}(G)$.

Firstly, according to Theorems 2.1 and 2.5, we get $\displaystyle \rho(A_{f}(P_{n}))\leq f(2, 2)\rho(A(P_{n}))=2f(2, 2)\cos\frac{\pi}{n+1}<2f(2, 2)$.

Next, let $T$ be a tree of order $n$ and $T\ncong P_{n}$. Then $\Delta(T)\geq3$. We prove $\displaystyle \rho(A_{f}(T))>\rho(A_{f}(P_{n}))$.

$\textbf{Case 1.}$ $\Delta(T)\geq4$.

If $\Delta(T)\geq4$, then by Theorems 2.1, 2.2 and 2.5, it follows that $\displaystyle \rho(A_{f}(T))\geq \rho(A_{f}(S_{5}))=2f(1,4)\geq 2f(2,2)>\rho(A_{f}(P_{n})).$

$\textbf{Case 2.}$ $\Delta(T)=3$.

$\textbf{Subcase 2.1.}$ There is a vertex $v_{i}$ in $T$ such that $d_{i}=3$ and every vertex $v_{j}\in N(v_{i})$ with $d_{j}\geq2$.

According to Theorems 2.1 and 2.2, we obtain $\displaystyle \rho(A_{f}(T))\geq\rho(A_{f}(T_{1}))$ (as shown in Fig. 3).
 \begin{figure}[H]
  \centering
    \includegraphics[width=3cm]{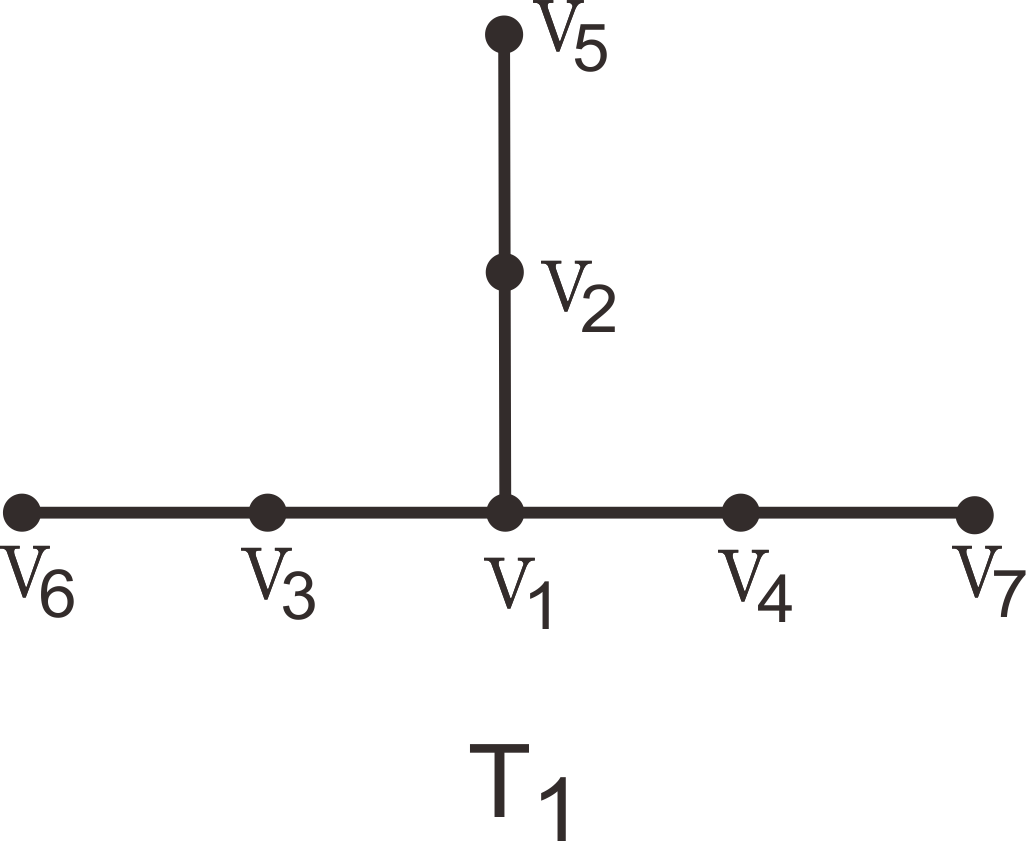}
  \caption*{\small \textbf{Fig. 3.} The tree $T_{1}$}
\end{figure}
The quotient matrix of the equitable partition $\{\{v_{1}\},\{v_{2},v_{3},v_{4}\},\{v_{5},v_{6},v_{7}\}\}$ of $A_{f}(T_{1})$ is
\begin{equation*}
  Q=\begin{pmatrix}
    0 & 3f(2,3) & 0\\
    f(2,3) & 0 & f(1,2) \\
 0 & f(1,2) & 0\\
\end{pmatrix}.
\end{equation*}

Let $\displaystyle P(x, T_{1})$ be the characteristic polynomial of the quotient matrix $Q$. By calculation, we get
\begin{center}
$\displaystyle P(x, T_{1})=x\big(x^{2}-3f^{2}(2, 3)-f^{2}(1, 2)\big)$.
\end{center}
Since $f(x, y)>0$ such that $\displaystyle f'_{x}(x, y)\geq0$ and $\displaystyle f''_{x}(x, y)\geq0$. Then $\displaystyle f^{2}(x, y)$ is convex in variable $x$ and we obtain $\displaystyle \sqrt{3f^{2}(2, 3)+f^{2}(1, 2)}\geq 2f(2, 2)$. It follows that $\displaystyle \rho(A_{f}(T))\geq \rho(A_{f}(T_{1}))\geq2f(2, 2)>\rho(A_{f}(P_{n}))$.

$\textbf{Subcase 2.2.}$ There is no vertex $v_{i}$ in $T$ such that $d_{i}=3$ and every vertex $v_{j}\in N(v_{i})$ with $d_{j}\geq2$.

Then $T$ is caterpillar tree with $\Delta(T)=3$. We can obtain $T$ by using Kelmans operation on the adjacent vertices of $P_{n}$ as shown in Fig. 4. Thus according to Theorem 3.1, we get $\displaystyle \rho(A_{f}(T))>\rho(A_{f}(P_{n}))$.

 \begin{figure}[H]
  \centering
    \includegraphics[width=9cm]{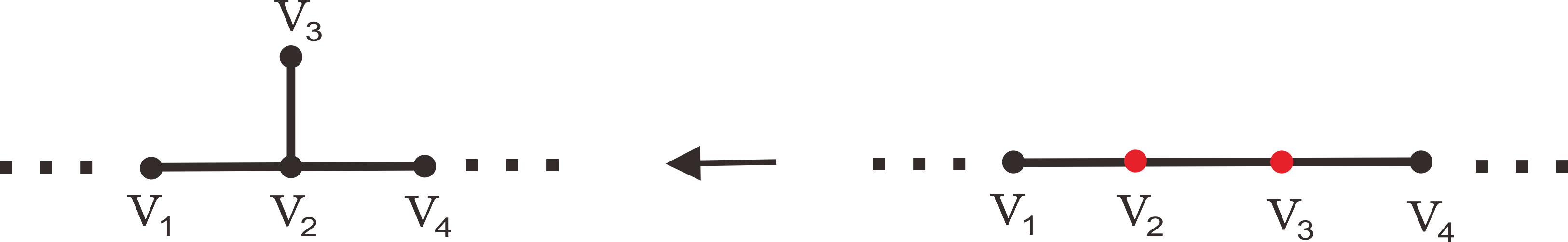}
      \caption*{\small \textbf{Fig. 4.} The Kelmans operation for subcase 2.2}
\end{figure}

So if $T\ncong P_{n}$, then $\displaystyle \rho(A_{f}(T))>\rho(A_{f}(P_{n}))$. This completes the proof of Theorem 4.1.
 \hfill$\Box$

\noindent $\textbf{Remark 4.1.}$ Li and Wang \cite{U9} tried to unify the spectral study of weighted adjacency matrices of graphs weighted by different topological indices. They found the extremal tree with the largest spectral radius of $A_{f}(G)$ is the star $S_{n}$ or the double star $S_{d,n-d}$ when $f(x, y)$ is increasing and convex in variable $x$. It is hard to determine which one is the largest one. We focus on restrictedly weighted adjacency matrices and push ahead Li and Wang's research.

\section{Conclusions and discussions}
\noindent

Restrictedly weighted adjacency matrices discussed in this paper include as special cases weighted adjacency matrices weighted by the following topological indices.

  1. First Zagreb index \cite{U17}: $\displaystyle f(x,y)=x+y$;

  2. First Hyper-Zagreb index \cite{U18}: $\displaystyle f(x,y)=(x+y)^{2}$;

  3. General Sum-Connectivity index \cite{U21}: $\displaystyle f(x,y)=(x+y)^{\alpha}$, where $\alpha\geq1$;

  4. Forgotten index \cite{U22}: $\displaystyle f(x,y)=x^{2}+y^{2}$;

  5. Somber index \cite{U4}: $\displaystyle f(x,y)=\sqrt{x^{2}+y^{2}}$;

  6. $p$-Sombor index \cite{U5}: $\displaystyle f(x,y)=\sqrt[p]{x^{p}+y^{p}}$ where $p\geq1$.

In this paper, we obtain $P_{n}$ is the unique tree with the smallest spectral radius, $S_{n}$ is the unique tree with the maximum spectral radius. If $f(x, y)$ is increasing and convex in variable $x$ but not possess the additional requirement in this paper, such as the following indices:

  1. Second Zagreb index \cite{U17}: $\displaystyle f(x,y)=xy$;

  2. Second hyper-Zagreb index \cite{10172}: $\displaystyle f(x,y)=(xy)^{2}$;

  3. First Gourava index \cite{10173}: $\displaystyle f(x,y)=x+y+xy$;

  4. Second Gourava index \cite{10173}: $\displaystyle f(x,y)=(x+y)xy$;

  5. First hyper-Gourava index \cite{10174}: $\displaystyle f(x,y)=(x+y+xy)^{2}$,

then Table 1 shows that we can not obtain a unified result.

\begin{table}[H]
   \caption{The approximate values of $\displaystyle \rho(A_{f}(G))$}
\renewcommand\arraystretch{2}
  \centering
  \begin{tabular}{ | c | l | l |l |l |l|l|l| }
    \hline
   $f(x,y)$  & $S_{15}$ & $S_{2,13}$ & $S_{3,12}$& $S_{4,11}$ & $S_{5,10}$ & $S_{6,9}$ & $S_{7,8}$\\
    \hline
  $xy$
     & 52.4
   & 52.0
   & 53.7
   & 56.4
   & 58.9
   & 60.9
   & \textbf{61.9}
    \\[3pt]
    \hline
    $(xy)^{2}$
   & 733.4
   & 894.3
   & 1381.3
   & 1973.6
   & 2518.4
   & 2926.1
   & \textbf{3142.9}
    \\[3pt]
    \hline
   $x+y+xy$
   & \textbf{108.5}
   & 102.1
   & 97.5
   & 94.2
   & 91.9
   & 90.5
   & 89.9
    \\[3pt]
    \hline
    $(x+y)xy$
   & 785.8
   & 741.4
   & 747.9
   & 781.5
   & 821.2
   & 853.8
   & \textbf{871.7}
    \\[3pt]
    \hline
    $(x+y+xy)^{2}$
   & 3146.7
   & 3033.7
   & 3326.4
   & 3864.2
   & 4433.3
   & 4883.3
   & \textbf{5127.7}
    \\[3pt]
    \hline
\end{tabular}
\end{table}

\section{Declaration of competing interest}

We declare that we have no conflicts of interest.

\section{Acknowledgements}
\noindent
The authors are grateful to Professor Xueliang Li for his talk in Xiamen which leads us to this research, and for some helpful comments on the first version of the paper. This work is supported by NSFC (No. 12171402).

\vskip0.5cm

\bibliographystyle{abbrv}
\bibliography{References}

\begin{thebibliography}{10}

\bibitem{U1}
A.~E. Brouwer and W.~H. Haemers.
\newblock Spectra of graphs.
\newblock 2011.

\bibitem{U8}
D.~K. Ch., G.~Ivan, M.~Igor, M.~Emina, and F.~Boris.
\newblock Degree–based energies of graphs.
\newblock {\em Lin. Algebra Appl.}, 554:185--204, 2018.

\bibitem{U13}
R.~Cruz, J.~Rada, and W.~Sanchez.
\newblock Extremal unicyclic graphs with respect to vertex-degree-based
  topological indices.
\newblock {\em MATCH Commun. Math. Comput. Chem.}, 88:481--503, 2022.

\bibitem{10}
E.~Estrada.
\newblock The {ABC} matrix.
\newblock {\em J. Math. Chem.}, 55(4):1021--1033, 2017.

\bibitem{8}
E.~Estrada, L.~Torres, L.~Rodríguez, and I.~Gutman.
\newblock An atom-bond connectivity index: Modelling the enthalpy of formation
  of alkanes.
\newblock {\em Indian J. Chem}, 37(10):849--855, 1998.

\bibitem{U22}
B.~Furtula and I.~Gutman.
\newblock A forgotten topological index.
\newblock {\em J. Math. Chem.}, 53:1184--1190, 2015.

\bibitem{U16}
W.~Gao.
\newblock Trees with maximum vertex-degree-based topological indices.
\newblock {\em MATCH Commun. Math. Comput. Chem.}, 88:535--552, 2022.

\bibitem{10172}
W.~Gao, M.~R. Farahani, M.~K. Siddiqui, and M.~K. Jamil.
\newblock On the first and second zagreb and first and second hyper-zagreb
  indices of carbon nanocones cnc$_{k}[n]$.
\newblock {\em J. Comput. Theor. Nanos.}, 13:7475--7482, 2016.

\bibitem{U4}
I.~Gutman.
\newblock Some basic properties of sombor indices.
\newblock {\em Open J. Discr. Appl. Math.}, 4:1--3, 2021.

\bibitem{U17}
I.~Gutman and N.~Trinajsti\'{c}.
\newblock Graph theory and molecular orbitals. total $\pi$-electron energy of
  alternant hydrocarbons.
\newblock {\em Chem. Phys. Lett.}, 17(4):535--538, 1972.

\bibitem{19}
R.~A. Horn and C.~R. Johnson.
\newblock {\em Matrix analysis}.
\newblock Cambridge University Press, Cambridge, 1990.

\bibitem{U14}
Z.~Hu, L.~Li, X.~Li, and D.~Peng.
\newblock Extremal graphs for topological index defined by a degree-based
  edge-weight function.
\newblock {\em MATCH Commun. Math. Comput. Chem.}, 88:505--520, 2022.

\bibitem{U15}
Z.~Hu, X.~Li, and D.~Peng.
\newblock Graphs with minimum vertex-degree function-index for convex
  functions.
\newblock {\em MATCH Commun. Math. Comput. Chem.}, 88:521--533, 2022.

\bibitem{U3}
A.~K. Kelmans.
\newblock On graphs with randomly deleted edges.
\newblock {\em Acta Math. Acad. Sci. Hung.}, 37:77--88, 1981.

\bibitem{10173}
V.~R. Kulli.
\newblock The gourava indices and coindices of graphs.
\newblock {\em Annals of Pure and Applied Mathematics}, 14(1):33--38, 2017.

\bibitem{10174}
V.~R. Kulli.
\newblock On hyper-gourava indices and coindices.
\newblock {\em International Journal of Mathematical Archive}, 8:116--120,
  2017.

\bibitem{U7}
X.~Li.
\newblock Indices, polynomials and matrices-a unified viewpoint.
\newblock {\em Invited talk at the $8$th Slovinian Conf, Graph Theory, Kranjska
  Gora}, pages June 21--27, 2015.

\bibitem{U10}
X.~Li, Y.~Li, and J.~Song.
\newblock The asymptotic value of graph energy for random graphs with
  degree-based weights.
\newblock {\em Discret. Appl. Math.}, 284:481--488, 2020.

\bibitem{U11}
X.~Li, Y.~Li, and Z.~Wang.
\newblock The asymptotic value of energy for matrices with
  degree-distance-based entries of random graphs.
\newblock {\em Lin. Algebra Appl.}, 603:390--401, 2020.

\bibitem{U12}
X.~Li, Y.~Li, and Z.~Wang.
\newblock Asymptotic values of four laplacian-type energies for matrices with
  degree distance-based entries of random graphs.
\newblock {\em Lin. Algebra Appl.}, 612:318--333, 2021.

\bibitem{10161}
X.~Li and D.~Peng.
\newblock Extremal problems for graphical function-indices and f-weighted
  adjacency matrix.
\newblock {\em Discrete Math. Lett.}, 9:57--66, 2022.

\bibitem{U9}
X.~Li and Z.~Wang.
\newblock Trees with extremal spectral radius of weighted adjacency matrices
  among trees weighted by degree-based indices.
\newblock {\em Lin. Algebra Appl.}, 620:61--75, 2021.

\bibitem{U6}
H.~Liu, L.~You, Y.~Huang, and X.~Fang.
\newblock Spectral properties of p-sombor matrices and beyond.
\newblock {\em MATCH Commun. Math. Comput. Chem.}, 87:59--87, 2021.

\bibitem{T1}
L.~M. Lov{\'a}sz and J.~Pelik{\'a}n.
\newblock On the eigenvalues of trees.
\newblock {\em Period. Math. Hungar.}, 3:175--182, 1973.

\bibitem{1}
M.~Randi\'{c}.
\newblock On characterization of molecular branching.
\newblock {\em J. Am. Chem. Soc.}, 97(23):6609--6615, 1975.

\bibitem{U5}
T.~R{\'e}ti, T.~Dosli{\'c}, and A.~Ali.
\newblock On the sombor index of graphs.
\newblock {\em Contrib. Math.}, 3:11--18, 2021.

\bibitem{555}
J.~A. Rodr\'{\i}guez.
\newblock A spectral approach to the {R}andi\'{c} index.
\newblock {\em Lin. Algebra. Appl}, 400:339--344, 2005.

\bibitem{666}
J.~A. Rodr\'{\i}guez and J.~M. Sigarreta.
\newblock On the {R}andi\'{c} index and conditional parameters of a graph.
\newblock {\em MATCH Commun. Math. Comput. Chem.}, 54(2):403--416, 2005.

\bibitem{15}
V.~S. Shegehall and R.~Kanabur.
\newblock Arithmetic-geometric indices of path graph.
\newblock {\em J. Math. Comput. Sci}, 6:19--24, 01 2015.

\bibitem{U18}
G.~H. Shirdel, H.~Rezapour, and A.~Sayadi.
\newblock The hyper-zagreb index of graph operations.
\newblock {\em Iran. J. Math. Chem.}, 4:213--220, 2013.

\bibitem{16}
L.~Zheng, G.~X. Tian, and S.~Y. Cui.
\newblock On spectral radius and energy of arithmetic-geometric matrix of
  graphs.
\newblock {\em MATCH Commun. Math. Comput. Chem.}, 83:635--650, 2020.

\bibitem{U21}
B.~Zhou and N.~Trinajstic.
\newblock On general sum-connectivity index.
\newblock {\em J. Math. Chem.}, 47:210--218, 2010.

\end{thebibliography}

\end{document}